\documentclass[oneside,reqno,english]{amsart}
\usepackage[T1]{fontenc}
\usepackage[utf8]{inputenc}
\usepackage{xcolor}
\usepackage{babel}
\usepackage{prettyref}
\usepackage{amstext}
\usepackage{amsthm}
\usepackage{amssymb}
\usepackage[pdfusetitle,
 bookmarks=true,bookmarksnumbered=false,bookmarksopen=false,
 breaklinks=false,pdfborder={0 0 0},pdfborderstyle={},backref=false,colorlinks=false]
 {hyperref}
\hypersetup{
 colorlinks=true,citecolor=blue,linkcolor=blue,linktocpage=true}

\makeatletter
\numberwithin{equation}{section}
\numberwithin{figure}{section}

\usepackage{prettyref}

\newrefformat{cor}{Corollary~\ref{#1}}
\newrefformat{subsec}{Section~\ref{#1}}
\newrefformat{lem}{Lemma~\ref{#1}}
\newrefformat{thm}{Theorem~\ref{#1}}
\newrefformat{sec}{Section~\ref{#1}}
\newrefformat{chap}{Chapter~\ref{#1}}
\newrefformat{prop}{Proposition~\ref{#1}}
\newrefformat{exa}{Example~\ref{#1}}
\newrefformat{tab}{Table~\ref{#1}}
\newrefformat{rem}{Remark~\ref{#1}}
\newrefformat{def}{Definition~\ref{#1}}
\newrefformat{fig}{Figure~\ref{#1}}
\newrefformat{claim}{Claim~\ref{#1}}
\newrefformat{assu}{Assumption~\ref{#1}}

\makeatother

\providecommand\theoremname{Theorem}
\theoremstyle{plain}
\newtheorem{thm}{\protect\theoremname}[section]
\providecommand\questionname{Question}
\newtheorem{question}[thm]{\protect\questionname}
\providecommand\lemmaname{Lemma}
\newtheorem{lem}[thm]{\protect\lemmaname}
\providecommand\definitionname{Definition}
\theoremstyle{definition}
\newtheorem{defn}[thm]{\protect\definitionname}
\providecommand\propositionname{Proposition}
\theoremstyle{plain}
\newtheorem{prop}[thm]{\protect\propositionname}
\providecommand\corollaryname{Corollary}
\newtheorem{cor}[thm]{\protect\corollaryname}
\begin{document}
\title{Joint Moments of Bilinear Resolvents}
\begin{abstract}
Given $A,B\in B\left(H\right)_{+}$, we consider the existence and
uniqueness of a positive operator $C\in B\left(H\right)_{+}$ such
that 
\[
\left(I+mA+nB+mnC\right)^{-1}=\int_{\left[0,1\right]^{2}}s^{m}t^{n}dE\left(s,t\right),\quad m,n\geq0,
\]
for a positive operator-valued measure $E$. This is a joint Hausdorff
moment problem. Without $C$, it holds if and only if $AB=BA$. In
dimension two we give a complete criterion for noncommuting pairs.
In higher finite dimensions, nonexistence is generic. Under certain
spectral and cyclicity assumptions, we obtain a necessary and sufficient
condition on arbitrary Hilbert spaces. The completion is unique in
this case. We also classify these pairs up to unitary equivalence.
The class includes irreducible infinite-dimensional examples, even
with continuous spectrum.
\end{abstract}

\author{James Tian}
\address{Mathematical Reviews, 535 W. William St, Suite 210, Ann Arbor, MI
48103, USA}
\email{james.ftian@gmail.com}

\maketitle
\tableofcontents{}

\section{Introduction}\label{sec:1}

The complete monotonicity of reciprocals of scalar polynomials has
been studied extensively \cite{MR0377423,MR1446489,MR3286037}. In
particular (see \cite[Theorem 3.1]{ACN2024}), for $a,b>0$ and $c\geq0$,
\begin{gather}
\left\{ \left(1+am+bn+cmn\right)^{-1}\right\} _{m,n\geq0}\;\text{has a joint Hausdorff representation}\nonumber \\
\Updownarrow\label{eq:1-1}\\
c\leq ab\nonumber 
\end{gather}

Here we fix $A,B\in B\left(H\right)_{+}$ and ask whether a positive
operator $C$ can be chosen so that the inverse has joint operator-valued
moments. The above scalar criterion \prettyref{eq:1-1} does not extend
directly. Specifically, we consider the following.
\begin{question}[positive completion]
Given $A,B\in B\left(H\right)_{+}$, can we find $C\in B\left(H\right)_{+}$
such that 
\begin{equation}
\left(I+mA+nB+mnC\right)^{-1}=\int_{\left[0,1\right]^{2}}s^{m}t^{n}dE\left(s,t\right),\quad m,n\geq0,\label{eq:1-2}
\end{equation}
for a positive operator-valued measure (POVM) $E$?
\end{question}

For $A,B\in B\left(H\right)_{+}$, without the mixed term, we show
in \prettyref{thm:2-5} that 
\begin{gather*}
\left\{ \left(I+mA+nB\right)^{-1}\right\} _{m,n\geq0}\text{ has a joint Hausdorff representation}\\
\Updownarrow\\
AB=BA
\end{gather*}
So the point of adding $C$ is to allow noncommuting pairs, although
its existence imposes much stronger spectral restrictions.

For a single operator $A\in B\left(H\right)_{+}$, the resolvents
$\left(I+mA\right)^{-1}$, $m\geq0$, form an operator-valued Hausdorff
moment sequence. The representing measure is determined by the spectral
measure of $A$. When \prettyref{eq:1-2} holds, the two one-variable
resolvents, $\left(I+mA\right)^{-1}$ and $\left(I+nB\right)^{-1}$,
are then its marginals, and the mixed coefficient $C$ accounts for
all the joint moments at once.

The question considered here is not about a simultaneous spectral
theorem. A POVM need not be projection-valued, and its marginals need
not commute. Nor does positivity of 
\[
\left(I+xA+yB+xyC\right)^{-1}
\]
at each point guarantee that its values are joint moments. The latter
requires positivity of all mixed finite differences, an infinite collection
of operator inequalities. The desired inverse polynomial leaves little
freedom to satisfy them.

In dimension two, that freedom has a simple description (\prettyref{thm:3-1}).
For noncommuting positive definite matrices $A$ and $B$, a positive
completion exists if and only if 
\[
\frac{A}{\left\Vert A\right\Vert }+\frac{B}{\left\Vert B\right\Vert }\geq I.
\]
When it exists it is unique. The condition is an ordinary operator
inequality, although it decides positivity of the entire two-variable
moment sequence. The situation appears to change sharply with the
dimension. For every matrix size at least three, the pairs with no
completion contain an open dense set (\prettyref{cor:4-2}).

The exceptional pairs in higher dimensions have a rigid spectral structure,
and that same structure gives a positive completion on an arbitrary
Hilbert space (\prettyref{sec:4}). 

We prove the following (\prettyref{thm:4-1}): Suppose the smallest
spectral values $a,b$ of $A,B$ are simple isolated eigenvalues,
with unit eigenvectors $q,p$, and $p$ is cyclic for $A$. Write
$R=pq^{*}/\left(q^{*}p\right)$. The condition is 
\[
\left(A-aI\right)\left(B-bI\right)=\gamma\left(I-R\right),\quad0<\gamma\le ab.
\]
It is both necessary and sufficient. The completion is determined
by $C=BA-\gamma R$.

In particular, there are irreducible infinite-dimensional pairs with
a positive completion, even when the operators have continuous spectrum.
Up to unitary equivalence, the pairs in this class are described by
three positive numbers and one scalar measure (\prettyref{cor:4-3}).

\subsection*{Related literature}

The underlying one-variable question belongs to the classical Hausdorff
moment problem, in which complete monotonicity of a sequence characterizes
measures on a compact interval \cite{MR1544592,MR0005923,MR0008438}.
The treatment of positive definite functions on semigroups provides
a useful general setting for moment and Laplace representations \cite{MR0747302}.
These results explain the individual resolvent measures, but they
do not decide whether the inverse of a prescribed noncommutative polynomial
has positive joint moments.

Multivariable moment problems require positivity on an algebra of
polynomials, rather than separate positivity in each coordinate. Existence
and representation on semialgebraic sets have been developed through
work on positive polynomials and moment functionals \cite{MR1092173,MR1254128,MR2383959,MR2589247}.
Positive moment matrices, flat extensions, and dimensional extensions
give other ways of approaching the multivariable problem \cite{MR1445490,MR1709313,MR1867381}.
These theories generally start with an unknown measure or an arbitrary
moment functional. Our moments are already fixed, apart from one coefficient
inside the inverse. This is what makes their positivity an operator
problem rather than an unrestricted moment problem.

Positive operator-valued measures appear in the dilation theory of
Naimark and the operator extensions of Halmos \cite{MR0002714,MR0044036}.
Stinespring's representation theorem and Arveson's work put positive
operator-valued measures in the setting of completely positive maps
\cite{MR0069403,MR0253059}. Choi's matrix criteria are central to
that setting \cite{MR0376726}. Andô's theorem provides commuting
dilations for pairs of commuting contractions \cite{MR0155193}. None
of these results by itself produces a bilinear inverse pencil with
given resolvent marginals. In particular, compatibility of positive
operator-valued marginals is a weaker requirement than the one imposed
here.

Moment positivity also has a long connection with subnormal operators
and commuting tuples. The classical subnormality criterion and its
operator-theoretic development may be found in \cite{MR0068129,MR1112128}.
The several-operator setting is represented, for example, by joint
isometries \cite{MR1066811}. Positive kernels give another language
for the same underlying positivity, starting with the reproducing
kernel theory of Aronszajn \cite{MR0051437} and its operator-valued
and several-variable developments \cite{MR3526117,MR1637941}. The
literature on multivariable operator monotone functions and Pick functions
is related through positivity of analytic operator functions \cite{MR2979860,MR3145735}.
These connections explain why joint moment positivity matters, but
they do not determine which resolvent pencils have it.

A different approach studies complete monotonicity of reciprocal scalar
polynomials. Early positive-coefficient questions include Szegő's
work \cite{MR1545428}. The theory of Riesz distributions on homogeneous
cones, and its treatment in harmonic analysis, gives striking criteria
for positivity of inverse powers of determinants \cite{MR0377423,MR1446489}.
Scott and Sokal developed this direction for several classes of polynomials
and related it to classical coefficient-positivity questions \cite{MR3286037}.
In the specific scalar bilinear case, a criterion for the reciprocal
to have joint Hausdorff moments is available \cite{ACN2024}. Our
problem differs at the first step. The coefficients are positive operators
which need not commute, and positivity of a scalar reciprocal, or
of the matrix entries of the inverse, does not give positivity of
an operator-valued representing measure. The spectral relation and
the operator-valued Laplace construction are needed to recover it. 

\section{Preliminaries: A scalar moment representation}\label{sec:scalar}

The scalar joint Hausdorff moment problem for reciprocals of bilinear
polynomials has been studied in the literature; see \cite[Theorem 3.1]{ACN2024}.
We need only the existence of a representing measure under a particular
condition on the coefficients. We give a direct proof.
\begin{lem}
\label{lem:scalar-moment} Let $a,b>0$ and $0\leq c\leq ab$. There
exists a probability measure $\mu$ on $\left[0,1\right]^{2}$ such
that 
\[
\frac{1}{1+am+bn+cmn}=\int_{\left[0,1\right]^{2}}s^{m}t^{n}d\mu\left(s,t\right),\quad m,n\geq0.
\]
\end{lem}

\begin{proof}
We consider three cases.

If $c=0$, then 
\[
\frac{1}{1+am+bn}=\int^{\infty}_{0}e^{-r}\left(e^{-ar}\right)^{m}\left(e^{-br}\right)^{n}dr.
\]
The pushforward of the probability measure $e^{-r}\,dr$ under the
map 
\[
r\longmapsto\left(e^{-ar},e^{-br}\right)
\]
is a representing probability measure on $\left[0,1\right]^{2}$.

If $c=ab$, then 
\[
\frac{1}{1+am+bn+abmn}=\frac{1}{1+am}\frac{1}{1+bn}.
\]
For $a>0$, 
\[
\frac{1}{1+am}=\int^{1}_{0}s^{m}\frac{1}{a}s^{1/a-1}ds.
\]
The analogous formula holds for $b$. Taking the product of these
two probability measures gives the required joint representation.

Now suppose $0<c<ab$, and let 
\[
\delta=ab-c>0.
\]
We have 
\[
\begin{aligned}1+am+bn+cmn & =c\left(m+\frac{b}{c}\right)\left(n+\frac{a}{c}\right)-\frac{\delta}{c}\\
 & =c\left(m+\frac{b}{c}\right)\left(n+\frac{a}{c}\right)\left[1-\frac{\delta}{c^{2}\left(m+b/c\right)\left(n+a/c\right)}\right].
\end{aligned}
\]
Since $m,n\geq0$, 
\[
0\leq\frac{\delta}{c^{2}\left(m+b/c\right)\left(n+a/c\right)}\leq\frac{\delta}{ab}<1.
\]
We can therefore expand the last factor as a geometric series: 
\[
\frac{1}{1+am+bn+cmn}=\frac{1}{c}\sum^{\infty}_{k=0}\left(\frac{\delta}{c^{2}}\right)^{k}\left(m+\frac{b}{c}\right)^{-k-1}\left(n+\frac{a}{c}\right)^{-k-1}.
\]

For every $\alpha>0$ and integer $k\geq0$, substitution $r=-\log s$
gives 
\[
\begin{aligned}\left(m+\alpha\right)^{-k-1} & =\frac{1}{k!}\int^{\infty}_{0}r^{k}e^{-\left(m+\alpha\right)r}dr=\frac{1}{k!}\int^{1}_{0}s^{m+\alpha-1}\left(-\log s\right)^{k}ds.\end{aligned}
\]
Thus each factor in the series has a positive one-variable Hausdorff
representation.

More explicitly, for each $k\geq0$, define a positive measure $\mu_{k}$
on $\left[0,1\right]^{2}$ by 
\[
\begin{aligned}d\mu_{k}\left(s,t\right)= & \frac{1}{c}\left(\frac{\delta}{c^{2}}\right)^{k}\cdot\frac{s^{b/c-1}\left(-\log s\right)^{k}}{k!}\frac{t^{a/c-1}\left(-\log t\right)^{k}}{k!}dsdt.\end{aligned}
\]
Its moments are 
\[
\int_{\left[0,1\right]^{2}}s^{m}t^{n}d\mu_{k}\left(s,t\right)=\frac{1}{c}\left(\frac{\delta}{c^{2}}\right)^{k}\left(m+\frac{b}{c}\right)^{-k-1}\left(n+\frac{a}{c}\right)^{-k-1}.
\]

In particular, taking $m=n=0$ gives 
\[
\mu_{k}\left(\left[0,1\right]^{2}\right)=\frac{1}{c}\left(\frac{\delta}{c^{2}}\right)^{k}\left(\frac{b}{c}\right)^{-k-1}\left(\frac{a}{c}\right)^{-k-1}=\frac{c}{ab}\left(\frac{\delta}{ab}\right)^{k}.
\]
Since $\delta<ab$, 
\[
\sum^{\infty}_{k=0}\mu_{k}\left(\left[0,1\right]^{2}\right)=\frac{c}{ab}\sum^{\infty}_{k=0}\left(\frac{\delta}{ab}\right)^{k}=\frac{c}{ab-\delta}=1.
\]

Consequently, 
\[
\mu=\sum^{\infty}_{k=0}\mu_{k}
\]
defines a probability measure on $\left[0,1\right]^{2}$. All the
measures are positive, so we may integrate term by term. Using the
geometric series above, we get 
\[
\begin{aligned}\int_{\left[0,1\right]^{2}}s^{m}t^{n}d\mu\left(s,t\right) & =\sum^{\infty}_{k=0}\frac{1}{c}\left(\frac{\delta}{c^{2}}\right)^{k}\left(m+\frac{b}{c}\right)^{-k-1}\left(n+\frac{a}{c}\right)^{-k-1}\\
 & =\frac{1}{1+am+bn+cmn}.
\end{aligned}
\]
\end{proof}

\section{Resolvent measures}\label{sec:2}

Let $H$ be a complex Hilbert space. For positive operators $A,B,C\in B\left(H\right)_{+}$,
let 
\[
F\left(x,y\right)=\left(I+xA+yB+xyC\right)^{-1},\quad x,y\geq0.
\]

\begin{defn}
We call $C$ a \textit{positive completion} of $A,B$ if there is
a positive $B\left(H\right)$-valued measure $E$ on $\left[0,1\right]^{2}$
such that 
\begin{equation}
F\left(m,n\right)=\int_{\left[0,1\right]^{2}}s^{m}t^{n}dE\left(s,t\right),\quad m,n\geq0.\label{eq:2-1}
\end{equation}
\end{defn}

The measure has total mass $I$. It is unique, since polynomials are
dense in $C(\left[0,1\right]^{2})$ and the scalar measures obtained
by taking matrix coefficients are determined by their polynomial integrals.
The integrals below can all be understood in this weak sense.

The one-variable resolvent has a simple representing measure, as follows. 
\begin{lem}
\label{lem:2-2}For $\lambda\geq0$, define the probability measure
$\rho_{\lambda}$ on $\left[0,1\right]$ by 
\begin{equation}
d\rho_{\lambda}\left(s\right)=\begin{cases}
\frac{1}{\lambda}s^{1/\lambda-1}ds, & \lambda>0,\\
d\delta_{1}\left(s\right), & \lambda=0.
\end{cases}\label{eq:2-2}
\end{equation}
Then, for any $T\in B\left(H\right)_{+}$ with spectral measure $E_{T}$,
we have
\begin{equation}
\left(I+xT\right)^{-1}=\int_{\left[0,1\right]}s^{x}d\nu_{T}\left(s\right),\quad x\geq0\label{eq:2-3}
\end{equation}
where 
\begin{equation}
\nu_{T}\left(V\right)=\int_{\sigma\left(T\right)}\rho_{\lambda}\left(V\right)dE_{T}\left(\lambda\right)\label{eq:2-4}
\end{equation}
 for all Borel $V\subset\left[0,1\right]$. In particular, $\nu_{T}$
is a normalized POVM.

Moreover, for $u>0$, set 
\[
g_{u}\left(\lambda\right)=\begin{cases}
e^{-u/\lambda}, & \lambda>0,\\
0, & \lambda=0.
\end{cases}
\]
Then 
\[
\nu_{T}\left(\left[0,e^{-u}\right]\right)=g_{u}\left(T\right).
\]
\end{lem}

\begin{proof}
Using \prettyref{eq:2-2}, direct integration gives
\[
\int_{\left[0,1\right]}s^{x}d\rho_{\lambda}\left(s\right)=\frac{1}{1+x\lambda},\quad x\geq0.
\]
By the spectral theorem, 
\begin{align*}
\left(I+xT\right)^{-1} & =\int_{\sigma\left(T\right)}\frac{1}{1+x\lambda}dE_{T}\left(\lambda\right)\\
 & =\int_{\sigma\left(T\right)}\left(\int_{\left[0,1\right]}s^{x}d\rho_{\lambda}\left(s\right)\right)dE_{T}\left(\lambda\right)=\int_{\left[0,1\right]}s^{x}d\nu_{T}\left(s\right),
\end{align*}
which gives \prettyref{eq:2-3}. 

Since each $\rho_{\lambda}$ is a probability measure, $\nu_{T}$
as in \prettyref{eq:2-4} is positive and has total mass $I$. Countable
additivity follows by monotone convergence. Uniqueness follows from
the one-variable Hausdorff moment theorem applied to matrix coefficients.

For $u>0$ and $\lambda>0$, we have 
\[
\begin{aligned}\rho_{\lambda}\left(\left[0,e^{-u}\right]\right) & =\frac{1}{\lambda}\int^{e^{-u}}_{0}s^{1/\lambda-1}ds=e^{-u/\lambda}.\end{aligned}
\]
For $\lambda=0$, the same interval has $\rho_{0}$-measure zero,
since $\rho_{0}=\delta_{1}$. Thus 
\[
\begin{aligned}\nu_{T}\left(\left[0,e^{-u}\right]\right) & =\int_{\sigma\left(T\right)}\rho_{\lambda}\left(\left[0,e^{-u}\right]\right)dE_{T}\left(\lambda\right)\\
 & =\int_{\sigma\left(T\right)}g_{u}\left(\lambda\right)dE_{T}\left(\lambda\right)=g_{u}\left(T\right).
\end{aligned}
\]
The function $g_{u}$ is continuous at zero since $e^{-u/\lambda}\to0$
as $\lambda\downarrow0$.
\end{proof}

Next, we extend \prettyref{eq:2-1} to real parameters.
\begin{prop}
\label{prop:2-3}Suppose $A,B,C\geq0$ and $C$ is a positive completion
of $A,B$. Let $E$ be the joint representing measure. Then
\begin{equation}
F\left(x,y\right)=\int_{\left[0,1\right]^{2}}s^{x}t^{y}dE\left(s,t\right),\quad x,y\geq0.\label{eq:2-5}
\end{equation}
\end{prop}

\begin{proof}
Fix $n\in\mathbb{Z}_{\geq0}$. For every Borel $V\subseteq\left[0,1\right]$,
the positive measure 
\[
V\longmapsto\int_{V\times\left[0,1\right]}t^{n}dE\left(s,t\right)
\]
has moments 
\[
\begin{aligned}\int_{\left[0,1\right]^{2}}s^{m}t^{n}dE\left(s,t\right) & =F\left(m,n\right)\\
 & =\left(I+nB+m\left(A+nC\right)\right)^{-1}\\
 & =\left(I+nB\right)^{-1/2}\left(I+mK_{n}\right)^{-1}\left(I+nB\right)^{-1/2},
\end{aligned}
\]
where 
\[
K_{n}=\left(I+nB\right)^{-1/2}\left(A+nC\right)\left(I+nB\right)^{-1/2}\geq0.
\]

By \prettyref{lem:2-2}, the positive measure 
\[
V\longmapsto\left(I+nB\right)^{-1/2}\nu_{K_{n}}\left(V\right)\left(I+nB\right)^{-1/2}
\]
has these same moments. Uniqueness of the one-variable Hausdorff representation
therefore gives 
\[
\int_{V\times\left[0,1\right]}t^{n}dE\left(s,t\right)=\left(I+nB\right)^{-1/2}\nu_{K_{n}}\left(V\right)\left(I+nB\right)^{-1/2}.
\]

We can now integrate $s^{x}$ for any real $x\geq0$. Again using
\prettyref{lem:2-2}, we obtain 
\[
\begin{aligned}\int_{\left[0,1\right]^{2}}s^{x}t^{n}dE\left(s,t\right) & =\left(I+nB\right)^{-1/2}\left(I+xK_{n}\right)^{-1}\left(I+nB\right)^{-1/2}\\
 & =\left(I+xA+nB+xnC\right)^{-1}=F\left(x,n\right).
\end{aligned}
\]

Now fix $x\geq0$. Applying the same argument in the second coordinate,
with $I+xA$ and $B+xC$ in place of $I+nB$ and $A+nC$, extends
the representation from nonnegative integers $n$ to every real $y\geq0$.
Thus \prettyref{eq:2-5} holds.
\end{proof}

The existence of a positive completion imposes a rigidity condition
on the spectral subspaces of $A$ and $B$.
\begin{thm}
\label{thm:2-4} Suppose $A,B,C\in B\left(H\right)_{+}$ and $C$
is a positive completion of $A,B$. Let $r>0$, $r\notin\sigma\left(B\right)$,
and write $P=E_{B}\left(\left[0,r\right]\right)$. Then, for all sufficiently
small $x>0$, 
\begin{equation}
\left(B+xC\right)\left(I+xA\right)^{-1}PH\subseteq PH.\label{eq:2-6}
\end{equation}
The same assertion holds with $A$ and $B$ interchanged.
\end{thm}

\begin{proof}
For $x\geq0$, let
\begin{equation}
T_{x}=\left(I+xA\right)^{-1/2}\left(B+xC\right)\left(I+xA\right)^{-1/2}.\label{eq:2-7}
\end{equation}
Then $T_{x}\geq0$ and 
\[
\left(B+xC\right)\left(I+xA\right)^{-1}=\left(I+xA\right)^{1/2}T_{x}\left(I+xA\right)^{-1/2}.
\]
It follows that \prettyref{eq:2-6} holds $\Longleftrightarrow$
\begin{equation}
T_{x}\left(I+xA\right)^{-1/2}PH\subseteq\left(I+xA\right)^{-1/2}PH.\label{eq:2-8}
\end{equation}

Since $T_{x}\to B$ in norm as $x\downarrow0$, and $PH=E_{B}\left(\left[0,r\right]\right)H$,
we will show that 
\begin{equation}
\left(I+xA\right)^{-1/2}PH=E_{T_{x}}\left(\left[0,r\right]\right)H\label{eq:2-9}
\end{equation}
for all sufficiently small $x>0$. Then \prettyref{eq:2-8} follows
from this and the spectral theorem. 

Let $E$ be the joint representing measure. We will use it to compare
the spectral subspaces of $T_{x}$ and $B$.

Fix $x\geq0$. For every Borel $V\subseteq\left[0,1\right]$, consider
the positive measure 
\[
V\longmapsto\int_{\left[0,1\right]\times V}s^{x}dE\left(s,t\right).
\]
By \prettyref{prop:2-3}, its moments are 
\[
\begin{aligned}\int_{\left[0,1\right]^{2}}s^{x}t^{n}dE\left(s,t\right) & =F\left(x,n\right)\\
 & =\left(I+xA+n\left(B+xC\right)\right)^{-1}\\
 & =\left(I+xA\right)^{-1/2}\left(I+nT_{x}\right)^{-1}\left(I+xA\right)^{-1/2},
\end{aligned}
\]
for every nonnegative integer $n$.

Since $T_{x}\geq0$, \prettyref{lem:2-2} gives 
\[
\left(I+nT_{x}\right)^{-1}=\int_{\left[0,1\right]}t^{n}d\nu_{T_{x}}\left(t\right).
\]
Therefore the positive measure 
\[
V\longmapsto\left(I+xA\right)^{-1/2}\nu_{T_{x}}\left(V\right)\left(I+xA\right)^{-1/2}
\]
has these same moments. Uniqueness of the one-variable Hausdorff representation
gives 
\begin{equation}
\int_{\left[0,1\right]\times V}s^{x}dE\left(s,t\right)=\left(I+xA\right)^{-1/2}\nu_{T_{x}}\left(V\right)\left(I+xA\right)^{-1/2}.\label{eq:2-10}
\end{equation}

We can now use positivity of the joint measure. Since $0\leq s^{x}\leq1$,
for every Borel set $V\subseteq\left[0,1\right]$ we have 
\[
0\leq\int_{\left[0,1\right]\times V}\left(1-s^{x}\right)dE\left(s,t\right).
\]
Equivalently, 
\[
\int_{\left[0,1\right]\times V}s^{x}dE\left(s,t\right)\leq\int_{\left[0,1\right]\times V}dE\left(s,t\right).
\]
Using \prettyref{eq:2-10} on both sides, and recalling that $T_{0}=B$,
this becomes 
\[
\left(I+xA\right)^{-1/2}\nu_{T_{x}}\left(V\right)\left(I+xA\right)^{-1/2}\leq\nu_{B}\left(V\right).
\]

Take $V=\left[0,e^{-u}\right]$, where $u>0$. By \prettyref{lem:2-2},
\[
\nu_{T_{x}}\left(\left[0,e^{-u}\right]\right)=g_{u}\left(T_{x}\right),\quad\nu_{B}\left(\left[0,e^{-u}\right]\right)=g_{u}\left(B\right).
\]
We therefore have 
\begin{equation}
\left(I+xA\right)^{-1/2}g_{u}\left(T_{x}\right)\left(I+xA\right)^{-1/2}\leq g_{u}\left(B\right),\quad x,u>0.\label{eq:2-11}
\end{equation}

We use this inequality to prove the inclusion in \prettyref{eq:2-9}.
Let $v\in PH$ and fix $q>r$. Since $P=E_{B}\left(\left[0,r\right]\right)$,
the spectrum of $B$ on $PH$ lies in $\left[0,r\right]$.

For $0\leq\lambda\leq r$, the definition of $g_{u}$ gives 
\[
0\leq g_{u}\left(\lambda\right)\leq e^{-u/r}.
\]
Hence 
\[
\left\langle v,g_{u}\left(B\right)v\right\rangle \leq e^{-u/r}\left\Vert v\right\Vert ^{2}.
\]
Combining this with \prettyref{eq:2-11}, we get 
\[
\begin{aligned} & \left\langle \left(I+xA\right)^{-1/2}v,g_{u}\left(T_{x}\right)\left(I+xA\right)^{-1/2}v\right\rangle \leq e^{-u/r}\left\Vert v\right\Vert ^{2}.\end{aligned}
\]

On the other hand, for every $\lambda\geq q$, 
\[
g_{u}\left(\lambda\right)=e^{-u/\lambda}\geq e^{-u/q}.
\]
The spectral theorem therefore gives 
\[
e^{-u/q}\left\Vert E_{T_{x}}\left(\left[q,\infty\right)\right)\left(I+xA\right)^{-1/2}v\right\Vert ^{2}\leq\left\langle \left(I+xA\right)^{-1/2}v,g_{u}\left(T_{x}\right)\left(I+xA\right)^{-1/2}v\right\rangle .
\]
Combining the last two inequalities, 
\[
\left\Vert E_{T_{x}}\left(\left[q,\infty\right)\right)\left(I+xA\right)^{-1/2}v\right\Vert ^{2}\leq e^{-u\left(1/r-1/q\right)}\left\Vert v\right\Vert ^{2}.
\]

Since $q>r$, we have $1/r-1/q>0$. The left-hand side is independent
of $u$, whereas the right-hand side tends to zero as $u\to\infty$.
Thus 
\[
E_{T_{x}}\left(\left[q,\infty\right)\right)\left(I+xA\right)^{-1/2}v=0.
\]
This holds for every $q>r$. Letting $q\downarrow r$ and using strong
convergence of spectral projections, we obtain 
\[
E_{T_{x}}\left(\left(r,\infty\right)\right)\left(I+xA\right)^{-1/2}v=0.
\]
Consequently, 
\[
\left(I+xA\right)^{-1/2}v\in E_{T_{x}}\left(\left[0,r\right]\right)H.
\]
Since $v\in PH$ was arbitrary, 
\[
\left(I+xA\right)^{-1/2}PH\subseteq E_{T_{x}}\left(\left[0,r\right]\right)H.
\]

It remains to prove equality for sufficiently small $x$. Since $T_{x}\to B$
in norm and $r\notin\sigma\left(B\right)$, the spectral projections
\[
Q_{x}=E_{T_{x}}\left(\left[0,r\right]\right)
\]
converge in norm to $P$.

Let $\Pi_{x}$ be the orthogonal projection onto $\left(I+xA\right)^{-1/2}PH$.
We also have $\Pi_{x}\to P$ in norm. To see this, the operator 
\[
P\left(I+xA\right)^{-1}P
\]
is positive and invertible on $PH$, and its restriction to $PH$
converges in norm to $I_{PH}$. The orthogonal projection onto $\left(I+xA\right)^{-1/2}PH$
is 
\[
\Pi_{x}=\left(I+xA\right)^{-1/2}P\cdot\left(\left.P\left(I+xA\right)^{-1}P\right|_{PH}\right)^{-1}\cdot P\left(I+xA\right)^{-1/2}.
\]
Each factor converges in norm to the corresponding factor at $x=0$.
Hence $\Pi_{x}\to P$.

The inclusion established above gives 
\[
\Pi_{x}\leq Q_{x}.
\]
Thus $Q_{x}-\Pi_{x}$ is an orthogonal projection, and 
\[
\left\Vert Q_{x}-\Pi_{x}\right\Vert \leq\left\Vert Q_{x}-P\right\Vert +\left\Vert \Pi_{x}-P\right\Vert \longrightarrow0.
\]
Since every nonzero orthogonal projection has norm $1$, we have $Q_{x}=\Pi_{x}$
for all sufficiently small $x>0$. Therefore 
\[
\left(I+xA\right)^{-1/2}PH=E_{T_{x}}\left(\left[0,r\right]\right)H.
\]
This is \prettyref{eq:2-9}, and hence \prettyref{eq:2-6} follows.

The assertion with $A$ and $B$ interchanged follows by symmetry.
\end{proof}

The affine case ($C=0$) has a stronger rigidity property. The following
result holds on an arbitrary Hilbert space.
\begin{thm}
\label{thm:2-5} Let $H$ be a complex Hilbert space and let $A,B\in B\left(H\right)_{+}$.
There is a positive $B\left(H\right)$-valued measure $E$ on $\left[0,1\right]^{2}$
such that 
\[
\left(I+mA+nB\right)^{-1}=\int_{\left[0,1\right]^{2}}s^{m}t^{n}dE\left(s,t\right)\quad\left(m,n\geq0\right)
\]
if and only if $AB=BA$. 
\end{thm}

\begin{proof}
Suppose the measure exists, and put 
\[
F\left(x,y\right)=\left(I+xA+yB\right)^{-1}.
\]
\prettyref{prop:2-3} shows that the measure represents $F\left(x,y\right)$
for all real $x,y\geq0$. 

For $x,y>0$, differentiation under the integral gives 
\begin{equation}
\left(-1\right)^{k+1}\frac{\partial^{k+1}F}{\partial x\partial y^{k}}\left(x,y\right)\geq0\quad\left(k\geq0\right).\label{eq:B-12}
\end{equation}
The logarithmic factors are bounded after multiplication by $s^{x}t^{y}$,
so the differentiation is justified. We may then let $x$ and $y$
tend to zero, since the derivatives of $F$ extend continuously to
the origin.

Now 
\[
-\frac{\partial F}{\partial x}\left(0,y\right)=\left(I+yB\right)^{-1}A\left(I+yB\right)^{-1}.
\]
Expanding near zero gives 
\begin{align*}
\left(I+yB\right)^{-1}A\left(I+yB\right)^{-1} & =\sum\nolimits_{j,\ell\geq0}(-y)^{j+\ell}B^{j}AB^{\ell}\\
 & =\sum\nolimits^{\infty}_{k=0}(-y)^{k}\Big(\underset{=:H_{k}}{\underbrace{\sum\nolimits^{k}_{j=0}B^{j}AB^{k-j}}}\Big).
\end{align*}
Using \prettyref{eq:B-12}, we get 
\begin{gather*}
\left(-1\right)^{k+1}\frac{\partial^{k+1}F}{\partial x\partial y^{k}}\left(0,0\right)=k!H_{k}\\
\Downarrow\\
H_{k}\geq0,\quad\forall k\geq0.
\end{gather*}

We show that these inequalities force $\left[A,B\right]=0$. To see
this, consider 
\[
0<a<b\leq c,\quad ac<b^{2}
\]
and set 
\[
P=E_{B}\left(\left[0,a\right]\right),\quad Q=E_{B}\left(\left[b,c\right]\right).
\]
Then we have 
\[
\left\Vert B^{j}P\right\Vert \leq a^{j},\quad\left\Vert B^{j}Q\right\Vert \leq c^{j}.
\]
Note that 
\begin{align*}
\left\Vert PH_{k}P\right\Vert  & \leq\sum^{k}_{j=0}\left\Vert PB^{j}AB^{k-j}P\right\Vert \\
 & \leq\sum^{k}_{j=0}\left\Vert PB^{j}\right\Vert \left\Vert A\right\Vert \left\Vert B^{k-j}P\right\Vert \\
 & \leq\left(k+1\right)a^{k}\left\Vert A\right\Vert .
\end{align*}
Similarly, 
\[
\left\Vert QH_{k}Q\right\Vert \leq\left(k+1\right)c^{k}\left\Vert A\right\Vert .
\]
Therefore, 
\begin{align*}
\left\Vert PH_{k}Q\right\Vert ^{2} & =\left\Vert \left(H^{1/2}_{k}P\right)^{*}\left(H^{1/2}_{k}Q\right)\right\Vert ^{2}\\
 & \leq\left\Vert \left(H^{1/2}_{k}P\right)^{*}\right\Vert ^{2}\left\Vert \left(H^{1/2}_{k}Q\right)\right\Vert ^{2}\\
 & =\left\Vert PH_{k}P\right\Vert \left\Vert QH_{k}Q\right\Vert \leq\left(k+1\right)^{2}\left(ac\right)^{k}\left\Vert A\right\Vert ^{2}.
\end{align*}
That is, 
\[
\left\Vert PH_{k}Q\right\Vert \leq\left(k+1\right)\left(ac\right)^{k/2}\left\Vert A\right\Vert .
\]

Observe that, for $k\geq1$, 
\[
H_{k}=\sum^{k}_{j=0}B^{j}AB^{k-j}=AB^{k}+BH_{k-1}.
\]
Since $Q=E_{B}\left(\left[b,c\right]\right)$ and $b>0$, the restriction
$B|_{QH}$ is invertible, with 
\[
\left\Vert \left(B|_{QH}\right)^{-k}\right\Vert \leq b^{-k}.
\]
It follows that 
\[
\begin{aligned}PH_{k}Q\left(B|_{QH}\right)^{-k} & =PAB^{k}Q\left(B|_{QH}\right)^{-k}+PBH_{k-1}Q\left(B|_{QH}\right)^{-k}\\
 & =PAQ+PBH_{k-1}Q\left(B|_{QH}\right)^{-k}.
\end{aligned}
\]
Therefore, 
\[
\begin{aligned}\left\Vert PAQ\right\Vert  & \leq\left\Vert PH_{k}Q\left(B|_{QH}\right)^{-k}\right\Vert +\left\Vert PBH_{k-1}Q\left(B|_{QH}\right)^{-k}\right\Vert \\
 & \leq b^{-k}\left(\left\Vert PH_{k}Q\right\Vert +\left\Vert PBH_{k-1}Q\right\Vert \right)\\
 & \leq b^{-k}\left(\left\Vert PH_{k}Q\right\Vert +a\left\Vert PH_{k-1}Q\right\Vert \right)\\
 & \leq\left\Vert A\right\Vert \left(\left(k+1\right)\left(\frac{\sqrt{ac}}{b}\right)^{k}+\frac{a}{b}k\left(\frac{\sqrt{ac}}{b}\right)^{k-1}\right).
\end{aligned}
\]
Since $ac<b^{2}$, the right-hand side tends to zero as $k\to\infty$.
Hence $PAQ=0$.

We have proved that 
\[
E_{B}\left(\left[0,a\right]\right)AE_{B}\left(\left[b,c\right]\right)=0
\]
whenever $0<a<b\leq c$ and $ac<b^{2}$.

For fixed $a>0$ and $b>a$, subdividing $\left[b,\left\Vert B\right\Vert \right]$
into sufficiently small spectral intervals therefore gives 
\[
E_{B}\left(\left[0,a\right]\right)AE_{B}\left(\left[b,\infty\right)\right)=0.
\]
Letting $b\downarrow a$, we obtain 
\[
E_{B}\left(\left[0,a\right]\right)AE_{B}\left(\left(a,\infty\right)\right)=0.
\]
Taking adjoints gives the reverse off-diagonal identity. Thus $A$
commutes with $E_{B}\left(\left[0,a\right]\right)$ for every $a>0$.
By the spectral theorem, $AB=BA$.

Conversely, suppose $AB=BA$. Then 
\[
\left(I+mA+nB\right)^{-1}=\int^{\infty}_{0}e^{-t}e^{-mtA}e^{-ntB}dt.
\]
For every $t\geq0$, the commuting positive contractions $e^{-tA}$
and $e^{-tB}$ have a joint spectral measure on $\left[0,1\right]^{2}$.
Integrating these measures against $e^{-t}dt$ gives the required
positive operator-valued measure. 
\end{proof}

The following commutativity criterion, implicit in the proof of \prettyref{thm:2-5},
may well be folklore. However, we have not found this exact formulation
in the literature. Since it may be of independent interest, we state
it separately.
\begin{prop}
\label{prop:2-6} Let $A,B\in B\left(H\right)_{+}$. For $k\geq0$,
set 
\[
H_{k}=\sum^{k}_{j=0}B^{j}AB^{k-j}.
\]
Then the following holds: 
\[
AB=BA\Longleftrightarrow H_{k}\geq0\quad\left(k\geq0\right).
\]
\end{prop}

\section{Completion in dimension two}\label{sec:3}

In dimension two, the positive completion problem has a fairly rigid
answer for noncommuting matrices. Although the problem asks for a
positive matrix $C$ such that the two-parameter family 
\[
\left(I+mA+nB+mnC\right)^{-1}
\]
admits a joint Hausdorff moment representation, existence is determined
by a single matrix inequality involving only $A$ and $B$. Moreover,
whenever a completion exists, it is unique and given by an explicit
formula. The necessity follows from the spectral invariance established
in \prettyref{sec:2}, and sufficiency is proved by reducing the matrix-valued
moment problem to three scalar moment representations.
\begin{thm}
\label{thm:3-1} Let $A,B\in M_{2}\left(\mathbb{C}\right)$ be (strictly)
positive definite and noncommuting. The pair $A,B$ has a positive
completion if and only if 
\begin{equation}
\frac{A}{\left\Vert A\right\Vert }+\frac{B}{\left\Vert B\right\Vert }\geq I.\label{eq:c-1}
\end{equation}
In that case the completion is unique, and it is given by
\begin{equation}
C_{0}=\left\Vert B\right\Vert A+\left\Vert A\right\Vert B-\left\Vert A\right\Vert \left\Vert B\right\Vert I.\label{eq:3-1}
\end{equation}
\end{thm}

\begin{proof}
``\textbf{$\Longrightarrow$}'' Suppose first that $C$ is a positive
completion. 

Write 
\[
a=\left\Vert A\right\Vert ,\quad b=\left\Vert B\right\Vert .
\]
Let $P$ be the spectral projection of $B$ for its smallest eigenvalue.
Choose $r$ strictly between the two eigenvalues of $B$. By \prettyref{thm:2-4},
for all sufficiently small $x>0$, 
\[
\left(B+xC\right)\left(I+xA\right)^{-1}PH\subseteq PH.
\]
We have the identity 
\[
\left(B+xC\right)\left(I+xA\right)^{-1}=B+x\left(C-BA\right)\left(I+xA\right)^{-1}.
\]
Since $PH$ reduces $B$, it follows that 
\[
\left(C-BA\right)\left(I+xA\right)^{-1}PH\subseteq PH.
\]
Expanding $\left(I+xA\right)^{-1}$ near zero and comparing coefficients
gives 
\[
\left(C-BA\right)A^{k}PH\subseteq PH,\quad\forall k\geq0.
\]

The line $PH$ is cyclic for $A$. Indeed, if it were an eigenspace
of $A$, then both matrices would be diagonal with respect to $PH\oplus\left(I-P\right)H$,
contrary to $AB\ne BA$. Since $H=PH+APH$, we obtain 
\[
\left(C-BA\right)H\subseteq PH.
\]
The operator $B$ acts as $bI$ on $\left(I-P\right)H$, so 
\[
\left(I-P\right)\left(C-bA\right)=0.
\]
The operator $C_{0}$ in \prettyref{eq:3-1} satisfies the same identity,
since $C_{0}-bA=a\left(B-bI\right)$. Hence 
\[
\left(C-C_{0}\right)H\subseteq PH.
\]

Interchange $A$ and $B$. If $Q$ is the spectral projection of $A$
for its smallest eigenvalue, the same argument gives 
\[
\left(C-C_{0}\right)H\subseteq QH.
\]

But $PH\cap QH=\left\{ 0\right\} $, otherwise they would be a common
eigenspace of $A$ and $B$, and the two matrices would commute. Thus
$C=C_{0}$. In particular, positivity of $C$ gives $C_{0}\geq0$,
which is \prettyref{eq:c-1}. \\

``$\Longleftarrow$'' Conversely, assume that 
\[
C_{0}=bA+aB-abI\geq0.
\]
where $a=\left\Vert A\right\Vert $ and $b=\left\Vert B\right\Vert $.
We need to show that $C_{0}$ is a positive completion of $A,B$. 

Let $d,e$ be the smaller eigenvalues of $A,B$, respectively. Then
$0<d<a$ and $0<e<b$.

The definition of $C_{0}$ gives 
\begin{equation}
I+xA+yB+xyC_{0}=\left(1+ax\right)\left(1+by\right)\Big[\underset{=:M}{\underbrace{I-\frac{x\left(aI-A\right)}{1+ax}-\frac{y\left(bI-B\right)}{1+by}}}\Big].\label{eq:c-3}
\end{equation}
Both $aI-A$ and $bI-B$ have rank one, with traces $a-d$ and $b-e$. 

Recall that, for rank-one $2\times2$ matrices $U,V$, 
\begin{equation}
\det\left(I-sU-tV\right)=\left(1-s\mathrm{tr}U\right)\left(1-t\mathrm{tr}V\right)-st\mathrm{tr}\left(UV\right).\label{eq:c-4}
\end{equation}
Indeed, the identities 
\begin{align*}
\det\left(I-T\right) & =1-\mathrm{tr}T+\det T\\
2\det T & =\left(\mathrm{tr}T\right)^{2}-\mathrm{tr}\left(T^{2}\right)
\end{align*}
give this formula \prettyref{eq:c-4}, with $T:=sU+tV$, using $\det U=\det V=0$.

Apply \prettyref{eq:c-4} to $M$ from \prettyref{eq:c-3} as follows.
Write 
\[
M=I-sU-tV
\]
where $U=aI-A$, $V=bI-B$, $s=\frac{x}{1+ax}$, $t=\frac{y}{1+by}$.
We get 
\begin{align*}
\det M & =\left(1-\frac{\left(a-d\right)x}{1+ax}\right)\left(1-\frac{\left(b-e\right)y}{1+by}\right)\\
 & \qquad-\frac{xy}{\left(1+ax\right)\left(1+by\right)}\mathrm{tr}\left(\left(aI-A\right)\left(bI-B\right)\right)\\
 & =\frac{\left(1+dx\right)\left(1+ey\right)-xy\mathrm{tr}\left(\left(aI-A\right)\left(bI-B\right)\right)}{\left(1+ax\right)\left(1+by\right)}.
\end{align*}
Substitute this into \prettyref{eq:c-3}, so
\begin{align}
\det\left(I+xA+yB+xyC_{0}\right) & =\left(1+ax\right)^{2}\left(1+by\right)^{2}\det M\nonumber \\
 & =\left(1+ax\right)\left(1+by\right)\left(1+dx+ey+cxy\right)\label{eq:c-5}
\end{align}
where 
\begin{equation}
c:=de-\mathrm{tr}\left(\left(aI-A\right)\left(bI-B\right)\right).\label{eq:3-2}
\end{equation}

In \prettyref{eq:c-5}, comparing the coefficients of $x^{2}y^{2}$,
we get 
\[
\det C_{0}=abc.
\]
Since $C_{0}\geq0$, we have $c\geq0$. Also, $aI-A\geq0$ and $bI-B\geq0$,
so 
\[
\mathrm{tr}\left(\left(aI-A\right)\left(bI-B\right)\right)=\mathrm{tr}\left(\left(bI-B\right)^{1/2}\left(aI-A\right)\left(bI-B\right)^{1/2}\right)\geq0.
\]
Back to \prettyref{eq:3-2}, it follows that 
\[
0\leq c\leq de.
\]

By \prettyref{lem:scalar-moment}, there is a positive measure $\mu$
on $\left[0,1\right]^{2}$ such that 
\begin{equation}
\frac{1}{1+dm+en+cmn}=\int_{\left[0,1\right]^{2}}s^{m}t^{n}d\mu\left(s,t\right),\quad m,n\geq0.\label{eq:3-3}
\end{equation}

We next calculate the matrix inverse
\[
F\left(x,y\right)=\left(I+xA+yB+xyC_{0}\right)^{-1}=\frac{1}{\left(1+ax\right)\left(1+by\right)}M^{-1},
\]
using \prettyref{eq:c-3}. 

Recall that, for any invertible $2\times2$ matrix $T$, 
\[
T^{-1}=\frac{\left(\mathrm{tr}T\right)I-T}{\det T}.
\]
For the matrix $M=I-sU-tV$ as above, we have 
\[
\text{tr}M=2-\left(a-d\right)s-\left(b-e\right)t.
\]
Thus, 
\[
\begin{aligned}\left(\text{tr}M\right)I-M & =\left(2-\left(a-d\right)s-\left(b-e\right)t\right)I-\left(I-sU-tV\right)\\
 & =I-s\left(A-dI\right)-t\left(B-eI\right).
\end{aligned}
\]
Also, from \prettyref{eq:c-5}, 
\[
\det M=\frac{1+dx+ey+cxy}{\left(1+ax\right)\left(1+by\right)}.
\]
Therefore, 
\begin{align}
F\left(x,y\right) & =\frac{1}{\left(1+ax\right)\left(1+by\right)}M^{-1}\nonumber \\
 & =\frac{1}{\left(1+ax\right)\left(1+by\right)}\frac{\left(\text{tr}M\right)I-M}{\det M}\nonumber \\
 & =\frac{1}{1+dx+ey+cxy}\left(I-\frac{x\left(A-dI\right)}{1+ax}-\frac{y\left(B-eI\right)}{1+by}\right)\nonumber \\
 & =\frac{1}{1+dx+ey+cxy}\left(\frac{\left(bd+ae\right)I-C_{0}}{ab}+\frac{A-dI}{a\left(1+ax\right)}+\frac{B-eI}{b\left(1+by\right)}\right)\label{eq:3-4}
\end{align}

All three constant matrices in \prettyref{eq:3-4} are positive semidefinite:
\[
K_{0}=\frac{\left(bd+ae\right)I-C_{0}}{ab},\quad K_{1}=\frac{A-dI}{a},\quad K_{2}=\frac{B-eI}{b}.
\]
The last two are immediate because $d,e$ are the smallest eigenvalues
of $A,B$. Since $C_{0}=bA+aB-abI$, we have 
\begin{align*}
\text{tr}C_{0} & =b\left(a+d\right)+a\left(b+e\right)-2ab=bd+ae.
\end{align*}
Recall that $C_{0}\geq0$, so 
\[
\left(bd+ae\right)I-C_{0}=\left(\text{tr}C_{0}\right)I-C_{0}\geq0.
\]

It remains to show that the scalar factors in \prettyref{eq:3-4}
have positive joint moment representations. Specifically, we now have
\begin{align*}
F\left(m,n\right) & =\frac{K_{0}}{1+dm+en+cmn}+\frac{K_{1}}{\left(1+dm+en+cmn\right)\left(1+am\right)}\\
 & \qquad+\frac{K_{2}}{\left(1+dm+en+cmn\right)\left(1+bn\right)}.
\end{align*}
The first is represented by \prettyref{eq:3-3}. 

The function $\left(1+am\right)^{-1}$ is a one-variable Hausdorff
moment sequence, since 
\[
\begin{aligned}\frac{1}{1+am} & =\frac{1}{a}\int^{1}_{0}u^{m+1/a-1}du=\int^{1}_{0}u^{m}d\rho\left(u\right),\end{aligned}
\]
where $\rho$ is the probability measure on $\left[0,1\right]$ with
density $a^{-1}u^{1/a-1}$.

Let $\mu$ be the positive measure in \prettyref{eq:3-3}. Multiplying
the two moment representations gives 
\[
\begin{aligned}\frac{1}{\left(1+dm+en+cmn\right)\left(1+am\right)} & =\int_{\left[0,1\right]^{2}}\int^{1}_{0}s^{m}t^{n}u^{m}d\rho\left(u\right)d\mu\left(s,t\right)\\
 & =\int_{\left[0,1\right]^{2}}\int^{1}_{0}\left(su\right)^{m}t^{n}d\rho\left(u\right)d\mu\left(s,t\right).
\end{aligned}
\]
Let $\nu_{1}$ be the pushforward of $\mu\otimes\rho$ under the map
\[
\left(s,t,u\right)\longmapsto\left(su,t\right).
\]
Since this map takes $\left[0,1\right]^{3}$ into $\left[0,1\right]^{2}$,
the measure $\nu_{1}$ is a positive measure on $\left[0,1\right]^{2}$,
and 
\[
\frac{1}{\left(1+dm+en+cmn\right)\left(1+am\right)}=\int_{\left[0,1\right]^{2}}s^{m}t^{n}d\nu_{1}\left(s,t\right).
\]

Similarly, we get a positive measure $\nu_{2}$ on $\left[0,1\right]^{2}$
such that 
\[
\frac{1}{\left(1+dm+en+cmn\right)\left(1+bn\right)}=\int_{\left[0,1\right]^{2}}s^{m}t^{n}d\nu_{2}\left(s,t\right).
\]

Now define the matrix-valued measure 
\[
E\left(\Omega\right)=K_{0}\mu\left(\Omega\right)+K_{1}\nu_{1}\left(\Omega\right)+K_{2}\nu_{2}\left(\Omega\right)
\]
for every Borel set $\Omega\subseteq\left[0,1\right]^{2}$. Since
$K_{0},K_{1},K_{2}\geq0$ and the three scalar measures are positive,
$E$ is a positive matrix-valued measure. By the three representations
above, 
\[
F\left(m,n\right)=\int_{\left[0,1\right]^{2}}s^{m}t^{n}\,dE\left(s,t\right),\quad m,n\geq0.
\]

Moreover, 
\[
E\left(\left[0,1\right]^{2}\right)=F\left(0,0\right)=I.
\]
Thus $E$ is a normalized positive operator-valued measure, and $C_{0}$
is a positive completion.
\end{proof}

\section{Higher dimensions}\label{sec:4}

The argument in \prettyref{thm:2-4} gives a restriction on $A$ and
$B$ when their lowest spectral subspaces are one-dimensional. In
higher dimensions this restriction is also sufficient. The proof works
for bounded operators on an arbitrary Hilbert space.
\begin{thm}
\label{thm:4-1} Let $H$ be a complex Hilbert space of dimension
at least two, and let $A,B\in B\left(H\right)$ be strictly positive.
Suppose that their smallest spectral values $a,b$ are simple isolated
eigenvalues, with unit eigenvectors $q,p$, respectively. Assume that
$p$ is cyclic for $A$. Let 
\[
R=\frac{1}{\left\langle q,p\right\rangle }\left|p\left\rangle \right\langle q\right|
\]

Then the pair $A,B$ has a positive  completion if and only if there
is a number $0<\gamma\le ab$ such that 
\begin{equation}
\left(A-aI\right)\left(B-bI\right)=\gamma\left(I-R\right).\label{eq:4-1}
\end{equation}
In that case the completion is unique and equals 
\begin{equation}
C=BA-\gamma R=AB+BA-aB-bA+\left(ab-\gamma\right)I.\label{eq:4-2}
\end{equation}
Equivalently, the spectral relation \prettyref{eq:4-1} can be written
\begin{equation}
\left(A-aI\right)^{1/2}\left(B-bI\right)\left(A-aI\right)^{1/2}=\gamma\left(I-\left|q\left\rangle \right\langle q\right|\right).\label{eq:4-3}
\end{equation}
\end{thm}

\begin{proof}
Note that since $p$ is cyclic for $A$ and $q$ is an eigenvector
of $A$, we have $\left\langle q,p\right\rangle \ne0$. Thus $R$
is well-defined. Moreover, $R^{2}=R$, with range $\mathbb{C}p$ and
kernel $q^{\perp}$.

\textbf{1.} Suppose first that $C$ is a positive completion. Choose
a number strictly between $b$ and the rest of $\sigma\left(B\right)$.
By \prettyref{thm:2-4}, with $P$ the projection onto $\mathbb{C}p$,
we have 
\[
\left(B+xC\right)\left(I+xA\right)^{-1}\mathbb{C}p\subseteq\mathbb{C}p
\]
for all sufficiently small $x>0$. Since 
\[
\left(B+xC\right)\left(I+xA\right)^{-1}=B+x\left(C-BA\right)\left(I+xA\right)^{-1},
\]
and $Bp=bp$, we have 
\[
\left(C-BA\right)\left(I+xA\right)^{-1}p\in\mathbb{C}p.
\]
Expand $\left(I+xA\right)^{-1}$ at zero. Comparing coefficients gives
\[
\left(C-BA\right)A^{k}p\in\mathbb{C}p\quad\left(k\ge0\right).
\]
By cyclicity of $p$, 
\[
\left(C-BA\right)H\subseteq\mathbb{C}p.
\]
That is, 
\[
C-BA=\left|p\left\rangle \right\langle w\right|
\]
for some $w\in H$; or, by taking adjoint, 
\begin{equation}
C-AB=\left|w\left\rangle \right\langle p\right|.\label{eq:5-4}
\end{equation}

We need only one more consequence of the measure argument. Apply \prettyref{thm:2-4}
with $A,B$ interchanged, using the isolated eigenvalue $a$ and its
eigenvector $q$. We obtain 
\[
\left(C-AB\right)q\in\mathbb{C}q.
\]
Since $\left\langle q,p\right\rangle \neq0$, using \prettyref{eq:5-4},
\[
\left(C-AB\right)q=w\left\langle p,q\right\rangle \in\mathbb{C}q
\]
and so $w\in\mathbb{C}q$. 

We have therefore proved 
\begin{equation}
C-BA=-\gamma R\label{eq:5-5}
\end{equation}
for some scalar $\gamma$. It is real, since $Rp=p$ and 
\[
\left\langle p,\left(C-BA\right)p\right\rangle =\left\langle p,Cp\right\rangle -b\left\langle p,Ap\right\rangle 
\]
is real. Notice that we did not assume $q$ cyclic for $B$.

Taking adjoint in \prettyref{eq:5-5} gives 
\begin{equation}
C-AB=-\gamma R^{*}.\label{eq:5-6}
\end{equation}
Combining \prettyref{eq:5-5}--\prettyref{eq:5-6}, we get
\[
\Big[\underset{S}{\underbrace{A-aI}},\underset{T}{\underbrace{B-bI}}\Big]=\left[A,B\right]=\gamma\left(R^{*}-R\right).
\]
Also, $RS=0$, because the range of $S$ is contained in $q^{\perp}$,
and $SR^{*}=0$, because the range of $R^{*}$ is $\mathbb{C}q$.

Let 
\[
K=ST-\gamma\left(I-R\right).
\]
Then 
\[
\begin{aligned}\left[S,K\right] & =S\left[S,T\right]+\gamma\left[S,R\right]\\
 & =\gamma S\left(R^{*}-R\right)+\gamma\left(SR-RS\right)\\
 & =\gamma\left(SR^{*}-RS\right)=0.
\end{aligned}
\]
Since $Tp=0$ and $Rp=p$, we also have $Kp=0$. As $K$ commutes
with $S$, 
\[
KS^{k}p=S^{k}Kp=0
\]
for every $k\geq0$. The vector $p$ is cyclic for $S$, so $K=0$,
i.e., 
\begin{equation}
ST=\gamma\left(I-R\right).\label{eq:5-7}
\end{equation}
This is precisely \prettyref{eq:4-1}. 

We next show that $\gamma>0$. Multiplying \prettyref{eq:5-7} on
the right by $S$ and using $RS=0$, we obtain 
\[
STS=\gamma S.
\]
For any nonzero $x\in q^{\perp}$, 
\[
\left\langle Sx,TSx\right\rangle =\gamma\left\langle x,Sx\right\rangle .
\]
Here $\left\langle x,Sx\right\rangle >0$. Also, $Sx$ is a nonzero
vector in $q^{\perp}$, while $\ker T=\mathbb{C}p$ and $\left\langle q,p\right\rangle \ne0$.
Hence $\left\langle Sx,TSx\right\rangle >0$, giving $\gamma>0$.

To obtain the upper bound for $\gamma$, observe that $Aq=aq$ implies
\[
\left\langle q,A^{-1}p\right\rangle =\frac{1}{a}\left\langle q,p\right\rangle .
\]
Therefore $RA^{-1}p=a^{-1}p$, and 
\[
CA^{-1}p=\left(BA-\gamma R\right)A^{-1}p=\left(b-\frac{\gamma}{a}\right)p.
\]
Since $C\geq0$, 
\[
0\leq\left\langle A^{-1}p,CA^{-1}p\right\rangle =\frac{ab-\gamma}{a}\left\langle A^{-1}p,p\right\rangle .
\]
The last inner product is strictly positive. Thus $\gamma\leq ab$.

Substituting 
\[
R=I-\gamma^{-1}\left(A-aI\right)\left(B-bI\right)
\]
into $C=BA-\gamma R$ gives 
\[
C=AB+BA-aB-bA+\left(ab-\gamma\right)I.
\]

It remains to verify the equivalent spectral relation. Since $RS^{1/2}=0$,
multiplying \prettyref{eq:5-7} on the right by $S^{1/2}$ gives 
\[
STS^{1/2}=\gamma S^{1/2}.
\]
The operator $S^{1/2}$ is boundedly invertible on $q^{\perp}$, because
$a$ is an isolated eigenvalue of $A$. Consequently, 
\begin{equation}
S^{1/2}TS^{1/2}=\gamma\left(I-\left|q\right\rangle \left\langle q\right|\right).\label{eq:5-8}
\end{equation}
That is, \prettyref{eq:4-3} holds.

Conversely, suppose that \prettyref{eq:5-8} holds. Since $S^{1/2}$
is invertible on $q^{\perp}$, the identity implies 
\[
STx=\gamma x
\]
for every $x\in q^{\perp}$. Also, $STp=0$, since $Tp=0$. As 
\[
H=q^{\perp}\oplus\mathbb{C}p,
\]
and $R$ is the idempotent with kernel $q^{\perp}$ and $Rp=p$, we
recover \prettyref{eq:5-7}.

\textbf{2.} We next prove sufficiency. Suppose that \prettyref{eq:4-1}
holds with $0<\gamma\leq ab$. Define 
\[
C=BA-\gamma R.
\]
Expanding \prettyref{eq:4-1} gives 
\[
C=AB+BA-aB-bA+\left(ab-\gamma\right)I,
\]
so $C$ is self-adjoint.

Decompose $H=\mathbb{C}p\oplus p^{\perp}$ and write 
\[
A=\begin{pmatrix}u & v^{*}\\
v & A_{1}
\end{pmatrix},\quad B=\begin{pmatrix}b & 0\\
0 & B_{1}
\end{pmatrix}.
\]
Put $D=A_{1}-aI$. Since $a$ is a simple isolated bottom eigenvalue
of $A$, there is $\varepsilon>0$ such that 
\[
A-aI\geq\varepsilon\left(I-\left|q\right\rangle \left\langle q\right|\right).
\]
For $z\in p^{\perp}$, 
\[
\begin{aligned}\left\langle z,Dz\right\rangle  & =\left\langle z,\left(A-aI\right)z\right\rangle \\
 & \geq\varepsilon\left(\left\Vert z\right\Vert ^{2}-\left|\left\langle q,z\right\rangle \right|^{2}\right)\geq\varepsilon\left|\left\langle q,p\right\rangle \right|^{2}\left\Vert z\right\Vert ^{2}.
\end{aligned}
\]
Thus $D$ is strictly positive and boundedly invertible.

The equation $\left(A-aI\right)q=0$, together with $\left\langle q,p\right\rangle \ne0$,
shows that $q$ is a nonzero scalar multiple of $\left(1,-D^{-1}v\right)$.
Hence 
\[
u=a+\left\langle v,D^{-1}v\right\rangle ,\quad R=\begin{pmatrix}1 & -v^{*}D^{-1}\\
0 & 0
\end{pmatrix}.
\]
Comparing the lower-right blocks in \prettyref{eq:4-1} gives 
\[
D\left(B_{1}-bI\right)=\gamma I,
\]
so $B_{1}=bI+\gamma D^{-1}$. We have obtained 
\begin{equation}
A=\begin{pmatrix}a+\left\langle v,D^{-1}v\right\rangle  & v^{*}\\
v & aI+D
\end{pmatrix},\quad B=\begin{pmatrix}b & 0\\
0 & bI+\gamma D^{-1}
\end{pmatrix}.\label{eq:4-5}
\end{equation}

Set 
\[
c=ab-\gamma,\quad h\left(x,y\right)=1+ax+by+cxy.
\]
We have $0\leq c<ab$. Define 
\[
\begin{aligned}K_{0} & =aA^{-1}+bB^{-1}-I,\\
K_{1} & =I-aA^{-1},\\
K_{2} & =I-bB^{-1}.
\end{aligned}
\]
Clearly, $K_{1},K_{2}\geq0$, because $A\geq aI$ and $B\geq bI$.
We show that $K_{0}\geq0$ as well.

From \prettyref{eq:4-5}, the operators $A_{1}$ and $B_{1}$ commute
and satisfy 
\[
\left(A_{1}-aI\right)\left(B_{1}-bI\right)=\gamma I.
\]
Expanding and multiplying by $A^{-1}_{1}B^{-1}_{1}$ gives 
\[
aA^{-1}_{1}+bB^{-1}_{1}-I=cA^{-1}_{1}B^{-1}_{1}.
\]
Consequently, 
\[
K_{0}=a\left(A^{-1}-\begin{pmatrix}0 & 0\\
0 & A^{-1}_{1}
\end{pmatrix}\right)+c\begin{pmatrix}0 & 0\\
0 & A^{-1}_{1}B^{-1}_{1}
\end{pmatrix}.
\]
The second term is positive because $A_{1}$ and $B_{1}$ are commuting
positive operators. For the first term, the block inverse formula
gives 
\[
A^{-1}-\begin{pmatrix}0 & 0\\
0 & A^{-1}_{1}
\end{pmatrix}=\frac{1}{u-\left\langle v,A^{-1}_{1}v\right\rangle }\begin{pmatrix}1\\
-A^{-1}_{1}v
\end{pmatrix}\begin{pmatrix}1 & -v^{*}A^{-1}_{1}\end{pmatrix}\geq0.
\]
Here the denominator is strictly positive, being the Schur complement
of $A_{1}$ in $A$. Thus $K_{0}\geq0$.

We have proved 
\[
K_{0},K_{1},K_{2}\geq0,\quad K_{0}+K_{1}+K_{2}=I.
\]
Moreover, $K_{1}$ commutes with $A$, and $K_{2}$ commutes with
$B$.

Now consider 
\[
G\left(x,y\right)=K_{0}+K_{1}\left(I+xA\right)^{-1}+K_{2}\left(I+yB\right)^{-1}.
\]
Since the coefficients are positive and sum to $I$, and both resolvents
are strictly positive, $G\left(x,y\right)$ is strictly positive and
boundedly invertible for every $x,y\geq0$.

The elementary resolvent identity gives 
\[
G\left(x,y\right)=\left(1+ax\right)\left(I+xA\right)^{-1}+\left(1+by\right)\left(I+yB\right)^{-1}-I.
\]

We claim that 
\[
G\left(x,y\right)\left(I+xA+yB+xyC\right)=h\left(x,y\right)I.
\]
To see this, recall that 
\[
C=BA-\gamma R=AB-\gamma R^{*},
\]
and that $AR^{*}=aR^{*}$ and $BR=bR$. Thus 
\[
\left(1+ax\right)\left(I+xA\right)^{-1}\left(I+xA+yB+xyC\right)=\left(1+ax\right)\left(I+yB\right)-\gamma xyR^{*},
\]
while 
\[
\left(1+by\right)\left(I+yB\right)^{-1}\left(I+xA+yB+xyC\right)=\left(1+by\right)\left(I+xA\right)-\gamma xyR.
\]
Combining these identities gives 
\[
\begin{aligned}G\left(x,y\right)\left(I+xA+yB+xyC\right) & =\left(1+ax\right)\left(I+yB\right)+\left(1+by\right)\left(I+xA\right)\\
 & \quad-\left(I+xA\right)\left(I+yB\right)-\gamma xyR\\
 & =\left(1+ax+by\right)I+xy\left(aB+bA-AB-\gamma R\right)\\
 & =h\left(x,y\right)I,
\end{aligned}
\]
where the last equality follows from \prettyref{eq:4-1}.

Since $G\left(x,y\right)>0$ and $h\left(x,y\right)>0$, we obtain
\[
I+xA+yB+xyC=h\left(x,y\right)G\left(x,y\right)^{-1}>0.
\]
Taking $x=y=t$, dividing by $t^{2}$, and letting $t\to\infty$ gives
$C\geq0$.

We have also proved the resolvent formula 
\begin{equation}
\begin{aligned}F\left(x,y\right) & =\left(I+xA+yB+xyC\right)^{-1}\\
 & =\frac{1}{h\left(x,y\right)}\left[K_{0}+K_{1}\left(I+xA\right)^{-1}+K_{2}\left(I+yB\right)^{-1}\right].
\end{aligned}
\label{eq:4-11}
\end{equation}

It remains to obtain a joint moment representation. By \prettyref{lem:scalar-moment},
the scalar sequence 
\[
\{h\left(m,n\right)^{-1}\}_{m,n\geq0}
\]
has a positive joint moment representation. By \prettyref{lem:2-2},
the resolvent sequences $\{\left(I+mA\right)^{-1}\}$ and $\{\left(I+nB\right)^{-1}\}$
have normalized positive operator-valued representing measures.

Since $K_{1}$ commutes with $A$ and $K_{2}$ commutes with $B$,
the sequences 
\[
K_{1}\left(I+mA\right)^{-1},\quad K_{2}\left(I+nB\right)^{-1}
\]
also have positive operator-valued representing measures. Combining
these with the scalar representation, as in the proof of \prettyref{thm:3-1},
gives a positive operator-valued measure representing $F\left(m,n\right)$.
It is normalized because $F\left(0,0\right)=I$.

Thus $C$ is a positive completion of $A,B$.
\end{proof}

The condition in \prettyref{thm:4-1} is restrictive, even though
there are examples in every dimension. First we retain the generic
nonexistence result.
\begin{cor}
\label{cor:4-2} For every $d\ge3$, the pairs of positive definite
Hermitian $d\times d$ matrices with no positive  completion contain
an open dense subset. The same is true for real symmetric positive
definite matrices.

For $d=2$, the noncommuting pairs having a positive  completion contain
a nonempty open subset. 
\end{cor}

\begin{proof}
For $d\ge3$, consider first the pairs whose smallest eigenvalues
are simple and for which the bottom eigenvector of $B$ is cyclic
for $A$ and the bottom eigenvector of $A$ is cyclic for $B$. These
conditions describe an open dense set. For density, perturb the two
matrices to have simple spectra and then make a small change of their
relative eigenbases. The required overlaps are nonzero except on sets
with empty interior, and the same argument works for real orthogonal
changes of basis.

For a pair in this set, let $a,q$ be the bottom eigenvalue and eigenvector
of $A$, and let $b,p$ be those of $B$. If a completion exists,
\prettyref{eq:4-3} implies that 
\[
\left(A-aI\right)^{1/2}\left(B-bI\right)\left(A-aI\right)^{1/2}
\]
is a scalar operator on $q^{\perp}$. This condition is closed relative
to the open set under consideration. The lowest eigenvalues and their
spectral projections depend continuously on the matrices.

To see that the scalar condition has empty interior, keep $A$, $p$,
and $b$ fixed and vary only the lower-right block of $B$ on $p^{\perp}$.
Small self-adjoint perturbations preserve positivity, the simple bottom
eigenvalue, and the cyclicity conditions. The restriction of $\left(A-aI\right)^{1/2}$
defines an invertible map 
\[
T:p^{\perp}\longrightarrow q^{\perp}.
\]
Indeed, its kernel is $p^{\perp}\cap\mathbb{C}q=\left\{ 0\right\} $,
and both spaces have dimension $d-1$. If $B_{1}$ is the lower-right
block of $B$, then the restriction of the preceding operator to $q^{\perp}$
is 
\[
T\left(B_{1}-bI\right)T^{*}.
\]
As $B_{1}$ varies, this can undergo arbitrary small self-adjoint
perturbations on $q^{\perp}$. When $d\ge3$, the space $q^{\perp}$
has dimension at least two. We can therefore perturb the restriction
away from scalar operators. The pairs violating the scalar condition
form an open dense subset, and none has a completion.

For $d=2$, choose noncommuting rank-one projections $P,Q$ and put
$A=I+\varepsilon P$, $B=I+\varepsilon Q$, with $0<\varepsilon<1$.
Both matrices have norm $1+\varepsilon$, and the matrix in \prettyref{eq:3-1}
is 
\[
\left(1+\varepsilon\right)\left[\left(1-\varepsilon\right)I+\varepsilon\left(P+Q\right)\right]>0.
\]
Thus these pairs have a completion by \prettyref{thm:3-1}. Strict
positivity and noncommutativity persist under small perturbations,
giving a nonempty open set. 
\end{proof}

We can also describe the operators in \prettyref{thm:4-1} up to unitary
equivalence. This uses only the spectral theorem and the cyclic vector
in its hypotheses.
\begin{cor}
\label{cor:4-3} The pairs in \prettyref{thm:4-1} which have a positive
completion are, up to unitary equivalence, the following pairs. Take
$a,b>0$, $0<\gamma\le ab$, and a nonzero finite positive measure
$\nu$ supported in a compact subset of $\left(0,\infty\right)$.
On $H=\mathbb{C}\oplus L^{2}\left(\nu\right)$ define 
\[
\begin{aligned}A\left(z,f\right) & =\left(\left(a+\int\frac{1}{t}d\nu\left(t\right)\right)z+\int f\left(t\right)d\nu\left(t\right),\quad z+\left(a+t\right)f\left(t\right)\right),\\
B\left(z,f\right) & =\left(bz,\ \left(b+\frac{\gamma}{t}\right)f\left(t\right)\right).
\end{aligned}
\]
In the second component of $A\left(z,f\right)$, the term $z$ means
the constant function with value $z$. Every such pair has a positive
completion. Two such pairs are unitarily equivalent if and only if
their numbers $a,b,\gamma$ and measures $\nu$ agree. 
\end{cor}

\begin{proof}
For a completed pair, use the block form \prettyref{eq:4-5}. The
vector $v$ is cyclic for $D$. Otherwise, the orthogonal complement
of its cyclic subspace would be a nonzero reducing subspace for $D$
orthogonal to $v$. This would give a nonzero reducing subspace for
$A$ orthogonal to $p$, contradicting cyclicity of $p$. Conversely,
if $v$ is cyclic for $D$, then $p=\left(1,0\right)$ is cyclic for
$A$. Indeed, a reducing subspace for $A$ orthogonal to $p$ would
lie in $\left\{ 0\right\} \oplus p^{\perp}$, would reduce $D$, and
would be orthogonal to $v$. There is no such nonzero subspace.

Let $E_{D}$ be the spectral measure of $D$ and set 
\[
\nu\left(V\right)=\left\langle v,E_{D}\left(V\right)v\right\rangle 
\]
for Borel $V\subset\left(0,\infty\right)$. The cyclic spectral theorem
identifies $p^{\perp}$ with $L^{2}\left(\nu\right)$, takes $D$
to multiplication by $t$, and takes $v$ to the constant function
$1$. Substitution in \prettyref{eq:4-5} gives the displayed formulas.
Conversely, for a nonzero measure with the stated support, multiplication
by $t$ has the cyclic vector $1$, because polynomials are dense
in $L^{2}\left(\nu\right)$. The operators satisfy \prettyref{eq:4-1}
and therefore have a positive completion.

The numbers $a,b$ are the bottom eigenvalues. Since $R$ is determined
by their one-dimensional eigenspaces and $I-R\ne0$, \prettyref{eq:4-1}
uniquely determines $\gamma$. The measure $\nu$ is determined by
the compression of $A-aI$ to $p^{\perp}$ and its coupling vector
$v$. These constructions respect unitary equivalence of the pair.
Thus equal data give unitarily equivalent pairs, and equivalent pairs
have equal data. 
\end{proof}

Here are two examples. The first settles a three-dimensional case
which is not covered by the two-dimensional theorem. 
\[
A=\begin{pmatrix}6 & 2 & 1\\
2 & 5 & 0\\
1 & 0 & 4
\end{pmatrix},\quad B=\begin{pmatrix}1 & 0 & 0\\
0 & 2 & 0\\
0 & 0 & 3
\end{pmatrix}.
\]
In \prettyref{eq:4-5}, $a=3$, $b=1$, $\gamma=2$, $D=diag\left(2,1\right)$,
and $v=\left(2,1\right)$. Its positive completion is 
\[
C=\begin{pmatrix}4 & 4 & 3\\
4 & 10 & 0\\
3 & 0 & 12
\end{pmatrix}.
\]
The pair generates $M_{3}\left(\mathbb{C}\right)$.

For an infinite-dimensional example, let $H=\mathbb{C}\oplus L^{2}\left(\left[1,2\right]\right)$,
let $D$ be multiplication by $t$, and let $v$ be the constant function
$1$. In \prettyref{eq:4-5} take $a=b=1$ and $\gamma=1/2$. The
theorem gives a positive completion. This pair is irreducible. To
see this, an operator commuting with $B$ preserves its simple eigenspace
$\mathbb{C}p$. On its orthogonal complement it commutes with multiplication
by $t$. Commutation with $A$ then requires that the lower block
send the constant function to a scalar multiple of itself. Since $1$
is cyclic for multiplication by $t$, the operator is scalar on the
lower block as well. Hence the joint commutant consists of scalars.
The multiplication operator has continuous spectrum.

There are also operator-algebraic restrictions on these pairs. They
hold for every pair in \prettyref{thm:4-1} in dimension at least
two, whether finite or infinite.
\begin{cor}
\label{cor:4-4} For a completed pair under the hypotheses of \prettyref{thm:4-1},
\[
rank\left(AB-BA\right)=2.
\]
If $H$ is infinite dimensional, then 
\[
C^{*}\left(A,B\right)=C^{*}\left(A\right)+\mathcal{K}\left(H\right)
\]
and 
\[
\sigma_{\mathrm{ess}}\left(B\right)=\left\{ b+\frac{\gamma}{t-a}:t\in\sigma_{\mathrm{ess}}\left(A\right)\right\} .
\]
\end{cor}

\begin{proof}
From \prettyref{eq:4-2} and self-adjointness of $C$, 
\[
AB-BA=\gamma\left(R^{*}-R\right).
\]
The vectors $p$ and $q$ are not proportional. Otherwise $p$ would
be an eigenvector of $A$, and cyclicity of $p$ would force $H=\mathbb{C}p$.
Since $q^{*}p\ne0$, the operator $R^{*}-R$ has rank two. This proves
the first assertion.

Let $G$ be zero on $\mathbb{C}q$ and equal to $\left(A-aI\right)^{-1}$
on $q^{\perp}$. Since $a$ is isolated, $G$ is bounded and belongs
to $C^{*}\left(A\right)$. Multiplying \prettyref{eq:4-1} on the
left by $G$ gives 
\[
\left(I-qq^{*}\right)\left(B-bI\right)=\gamma G\left(I-R\right).
\]
Consequently, 
\[
B-\left(bI+\gamma G\right)=qq^{*}\left(B-bI\right)-\gamma GR.
\]
The right side has rank at most two. Hence $B$ is a finite-rank perturbation
of $bI+\gamma G$ and the essential-spectrum identity follows.

The pair acts irreducibly. Indeed, any operator commuting with both
$A$ and $B$ preserves the one-dimensional eigenspace $\mathbb{C}p$
and acts on it as a scalar. Since $p$ is cyclic for $A$, the operator
is scalar on $H$. Thus $C^{*}\left(A,B\right)$ is an irreducible
$C^{*}$-algebra. It contains the nonzero finite-rank commutator $AB-BA$.
The standard theorem on irreducible $C^{*}$-algebras containing a
nonzero compact operator gives $\mathcal{K}\left(H\right)\subseteq C^{*}\left(A,B\right)$.
The preceding finite-rank identity also gives $B\in C^{*}\left(A\right)+\mathcal{K}\left(H\right)$.
Combining the two inclusions proves the asserted equality. 
\end{proof}

\section{Final remark}

We end with two questions, one concerns the spectral structure of
pairs admitting positive completions, the other is about the dilation
of their representing measures.
\begin{enumerate}
\item \prettyref{thm:4-1} assumes that the smallest spectral values of
$A$ and $B$ are simple isolated eigenvalues. It is not clear whether
an irreducible pair of strictly positive operators can admit a positive
completion when neither operator has eigenvalues.
\item A positive completion also gives a simultaneous dilation. Indeed,
let $E$ be the representing POVM, then there exist a spectral measure
$P$ on $\left[0,1\right]^{2}$ and an isometry $V:H\longrightarrow K$
such that 
\[
E\left(U\right)=V^{*}P\left(U\right)V.
\]
Set 
\[
S=\int_{\left[0,1\right]^{2}}sdP\left(s,t\right),\quad T=\int_{\left[0,1\right]^{2}}tdP\left(s,t\right).
\]
Then $S,T$ are commuting positive contractions and 
\[
\left(I+mA+nB+mnC\right)^{-1}=V^{*}S^{m}T^{n}V,\quad m,n\geq0.
\]
We may take 
\[
K=\overline{span}\left\{ S^{m}T^{n}Vh:m,n\geq0,\ h\in H\right\} .
\]
This minimal dilation is unique up to unitary equivalence. It would
be interesting to determine its spectral multiplicity and the structure
of its spectral measure, particularly how these depend on the completion
$C$.
\end{enumerate}
\bibliographystyle{amsalpha}
\bibliography{ref}

\end{document}